\documentclass{article}
\usepackage{myfig,latexsym,amsthm,amsfonts,epsfig,psfrag}
\begin{document}

\newtheorem{lemma}{Lemma}
\newtheorem{theorem}{Theorem}

\def \H{{I\!\!H}}
\def \R{{I\!\!R}}
\def \Z{{\mathbb Z}}
\def \g{{\mathfrak g}}
\def \h{{\mathfrak h}}

\begin{center}
{\LARGE  \bf
Maximal rank root subsystems \\
\vspace{1.5pt}
of hyperbolic root systems.
}

\vspace{17pt}
{\large P.~Tumarkin}
\end{center}

\vspace{10pt}

\begin{center}
\parbox{10cm}{\scriptsize
{\it Abstract.} A Kac-Moody algebra is called 
hyperbolic  if it corresponds to 
a generalized Cartan matrix of hyperbolic type. We study root
subsystems of root systems of hyperbolic algebras. 
In this paper, we classify  maximal rank regular hyperbolic 
subalgebras of hyperbolic Kac-Moody algebras.}
\end{center}

\vspace{10pt}

\section*{Introduction}

A  generalized Cartan matrix $A$ is called a matrix of 
{\it hyperbolic type} if it is indecomposable symmetrizable
of indefinite type, and if any proper principal submatrix of the
corresponding symmetric matrix $B$ is of finite or affine type.
In this case $B$ is of the signature  $(n,1)$.

Consider a generalized Cartan matrix $A$ of  hyperbolic type.
Following Kac~\cite{Kac}, we can construct a Kac-Moody algebra
$\g(A)$.  According to Vinberg~\cite{Vinberg}, 
the Weyl group of the root system $\Delta(A)$ is a Coxeter group.
A fundamental chamber of the Weyl group is an $n$-dimensional
hyperbolic Coxeter simplex of finite volume,
whose dihedral angles are in the set $\{\frac{\pi}{2},
 \frac{\pi}{3}, \frac{\pi}{4},\frac{\pi}{6}\}$
(zero angle can also appear if $n=2$).

In analogy with the finite-dimensional theory (see~\cite{Dyn}), 
we say a subalgebra $\g_1\subset\g(A)$ to be  
{\it regular} if $g_1$ is invariant with respect to some Cartan
subalgebra $\h$ of  $\g(A)$.
In other words, $\g_1\subset\g(A)$ is regular if it has a basis
composed of some elements of $\h$  and some root vectors of  $\g(A)$
(with respect to  $\h$). 
We are interested in maximal rank regular subalgebras that can be
constructed as  Kac-Moody algebras $\g_1(A_1)$ for some 
 generalized Cartan matrix $A_1$ of hyperbolic type.

Any subalgebra of this type of the  Kac-Moody algebra $\g(A)$
has a root system
$\Delta_1 (A_1)\subset\Delta (A)$ such that
$$
{\mbox if }\ \alpha, \beta\in\Delta_1 \ {\mbox and \ }
\ \alpha+\beta\in\Delta,\
{\mbox then \ }\ \alpha+\beta\in\Delta_1\leqno{(*)}
$$
Conversely, suppose we have a hyperbolic root system $\Delta_1$  
in a hyperbolic root system $\Delta (A)$, and
($*$) holds. Then we can construct a subalgebra of  $\g(A)$ we are
interested in.

By 
{\it hyperbolic root system} we mean a root system of a Kac-Moody
algebra constructed on a generalized Cartan matrix of hyperbolic type.

Let $\Delta$ be a hyperbolic root system.
A root system $\Delta_1\subset\Delta $ is called 
a {\it root subsystem }  of $\Delta$ 
if the condition ($*$) holds.

$\!$The classification of root subsystems of finite root systems
is due to Dynkin$\!$~\cite{Dyn}.

In this paper we classify maximal rank hyperbolic root subsystems
of hyperbolic root systems.

Consider a  maximal rank hyperbolic root subsystem $\Delta_1$ 
of a hyperbolic root system $\Delta$.
Let  $W_1$ and $W$ be the Weyl groups of  $\Delta_1$ and $\Delta$ respectively.
Let  ${F_1}$ and $F$  be  fundamental chambers of $W_1$ and $W$.
Then  ${F_1}$ and $F$ are hyperbolic Coxeter simplices of finite volume.
The groups $W_1$ and $W$ are generated by the reflections with respect 
to the facets of  ${F_1}$ and $F$  respectively.
Since $W_1$ is a subgroup of $W$, the simplex ${F_1}$
is composed of several copies of $F$.
Moreover, any two copies of $F$ having a common facet
are symmetric with respect to this facet.

By {\it reflection group} we mean a group generated by reflections.
Introduce a partial ordering $\ge$ on the set of reflection 
subgroups of $W$ by setting $G\ge H$ if $H\subset G$.
A decomposition ($F,{F_1}$) of a simplex $F_1$ into several copies of
$F$ is called {\it minimal} if $W_1$ is a maximal proper reflection 
subgroup of $W$.
All the minimal decompositions of hyperbolic Coxeter simplices
of finite volume are listed in~\cite{Fel4-9},~\cite{Fel3} and~\cite{Klim}.

From now on by simplex we mean a hyperbolic Coxeter simplex of finite
volume, whose dihedral angles are in the set $\{\frac{\pi}{2},
 \frac{\pi}{3}, \frac{\pi}{4},\frac{\pi}{6}\}$
(zero angle can also appear if $n=2$). 

In Section~\ref{sublattice-max} (Th.~\ref{max}) we prove that 
any minimal decomposition of a hyperbolic simplex corresponds to some
root subsystem of a hyperbolic root system.  
In Section~\ref{sublattice} (Th.~\ref{real}) we prove that 
any decomposition of a hyperbolic simplex corresponds to some
root subsystem.  
The complete classification of maximal rank hyperbolic root subsystems
is contained in Fig.~\ref{2_1}--\ref{9}. 

The author is grateful to Prof. E.~B.~Vinberg for his attention to the
work and useful remarks.

\vspace{15pt}

\section{Maximal subgroups}\label{sublattice-max}

We use the following notation: 
$A$ is a generalized Cartan matrix of hyperbolic type;
$\Delta$ is the corresponding root system;
$\alpha_1,...,\alpha_{n+1}$ are simple roots;
$F$ is a fundamental chamber of $\Delta$; 
$L=\sum\limits_{i=1}^{n+1}\Z\alpha_i$ is the corresponding root lattice.
The Weyl group $W_F$ of $\Delta$ is generated by the reflections with
respect to the facets of the simplex $F$.
The simple roots vanish on the facets of $F$.
Furthermore, $\Delta^{\vee}$ and $L^{\vee}$ are the root system and
the root lattice for the generalized Cartan matrix $A^t$ 
(the fundamental simplex of the Weyl group of $\Delta^{\vee}$ 
is the same as of $\Delta$, but the lengths of simple roots are 
different in these systems). 
$\Delta_1\subset\Delta$ is a hyperbolic root system 
whose root lattice
$L_1$ is a maximal rank sublattice of $L$;  
${F_1}$ is a fundamental simplex of
the Weyl group $W_{F_1}$ 
of the root system $\Delta_1$.

We will use the following description of root system
(see~\cite{Kac}). 
A hyperbolic root system $\Delta$ consists of two disjoint parts: 
the set of real roots
$\Delta^{re}$ and the set of imaginary roots $\Delta^{im}$, where

$\Delta^{re}=W(\alpha_1)\bigcup\dots\bigcup W(\alpha_{n+1})\,$,\qquad
$\Delta^{im}=\{\alpha\in L\ |\ (\alpha|\alpha)\le 0\}$.

\begin{lemma}\label{equiv}
Let $(F,{F_1})$ be a minimal decomposition.
The following four conditions are equivalent:

\begin{itemize}
\item[{\sc (i)}] $L_1$ is a proper sublattice of $L$.

\item[{\sc (ii)}] $\Delta_1=\Delta\cap L_1$.

\item[{\sc (iii)}]  $\Delta_1$ is a root subsystem of $\Delta$.

\item[{\sc (iv)}] The condition {\sc ($*$)} holds for
$\Delta_1^{re}\subset\Delta^{re}$.

\end{itemize}
\end{lemma}

\begin{proof}\quad
{\sc (i)$\to$(ii)} \quad
For $\Delta_1^{im}$ the statement is evident.

Suppose that there exists $\alpha\in\Delta^{re}$ such that $\alpha\in
L$ and $\alpha\notin\Delta_1$. 
Consider a subgroup $G$ of $W_F$ generated by the reflections with
respect to all the roots contained in $L_1$. Clearly,
$W_{F_1}\subset G$. Since $W_{F_1}$ is maximal in
$W_F$ and  $G\ne W_{F_1}$, we have $G=W_F$. 
Hence, any simple root of $\Delta$ can be written as 
$\sum\limits_{i=1}^{n+1}c_i \beta_i$, where $\beta_i \in L_1$ and 
$c_i\in \Z$. Thus, any simple root of  $\Delta$ belongs to $L_1$,
and $L=L_1$.

\vspace{3pt}

{\sc (ii)$\to$(iii)} \quad 
Suppose that $\alpha, \beta\in\Delta_1$. 
Then $\alpha,\beta,\alpha+\beta\in L_1$.
If $\alpha+\beta\in\Delta$ then $\alpha+\beta\in\Delta\cap L_1$.
Therefore,
$\alpha+\beta\in\Delta_1$.

\vspace{3pt}

{\sc (iii)$\to$(iv)} \quad
$\!\!$The proof is evident.

\vspace{3pt}

{\sc (iv)$\to$(i)} \quad
$\;$Assume that $L=L_1$.

Suppose that the simplex $F_1$ has a decomposed dihedral angle,
i.e. some mirror of a reflection contained in $W_F$
decomposes the dihedral angle of $F_1$.
Let $\alpha$ and $\beta\in\Delta_1^{re}$ be the roots vanishing on the
facets of this dihedral angle (the roots are the outward normals to
the facets of the angle). Then $\alpha-\beta\in L_1$ 
vanishes on one of the mirrors decomposing the dihedral angle.
Hence, there exists $c>0$ such that $c(\alpha-\beta)\in\Delta^{re}$.
By the assumption $L=L_1$, thus, 
$c(\alpha-\beta)\in L_1$.
Without loss of generality we can assume that
$\alpha$ and $\beta$ are simple roots of  $\Delta_1$.  
Then $\alpha-\beta\notin\Delta_1^{re}$.
The lattice $L_1$ is generated by simple roots of
$\Delta_1$, thus, $c\in\Z$.
Since  ($*$) holds for $\Delta_1^{re}\subset\Delta^{re}$ 
and $\alpha-\beta\notin\Delta_1^{re}$, we have $c\ne 1$. 
Since $(\alpha\ |\ \beta)\le 0$, if $c\ge 2$ then
$c(\alpha-\beta)$ is more than two times longer than $\alpha$.
This is impossible, since $c(\alpha-\beta)$ and  $\alpha$ are not
mutually orthogonal.

Suppose now that ${F_1}$ has no decomposed dihedral angle.
Then the pair $(F,{F_1})$ is one of the six pairs listed in 
Table~\ref{fund}.
Since no of these simplices has a dihedral angle different from
$\frac{\pi}{2}$ and $\frac{\pi}{3}$,
each simplex contained in Table~\ref{fund} corresponds to a unique
root system. A direct calculation shows that for each of these six pairs
 roots of the subsystem generate an index two sublattice of
the root lattice.

\end{proof}

\begin{table}[htb!]
\begin{center}
\begin{tabular}{|c|c|c|}
\hline
$\vphantom{\int_a^A}F\vphantom{\int_a^A}$& ${F_1}$ & $[W_F:W_{F_1}]$  \\
\hline
&&\\
\raisebox{2pt}{
\epsfig{file=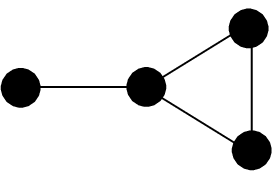,width=0.065\linewidth}}&
{\epsfig{file=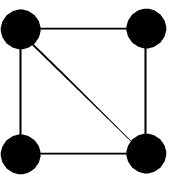,width=0.055\linewidth}}&
\raisebox{5pt}{5}
 \\
\raisebox{1pt}{
\epsfig{file=f3-1.eps,width=0.065\linewidth}}&
\raisebox{-2pt}{\epsfig{file=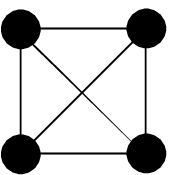,width=0.06\linewidth}}&
\raisebox{4.5pt}{12} 
\\
\raisebox{1pt}{
\epsfig{file=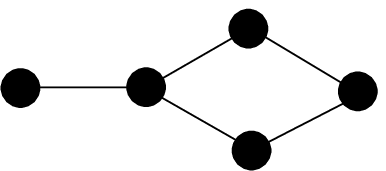,width=0.08\linewidth}}&
$\vphantom{\int\limits_a^{A^A}}$\epsfig{file=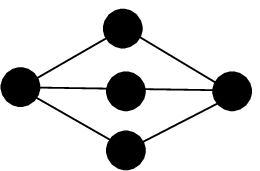,width=0.06\linewidth}$\vphantom{\int\limits_a^{A^A}}$&
\raisebox{3pt}{10} 
\\
\raisebox{1pt}{
\epsfig{file=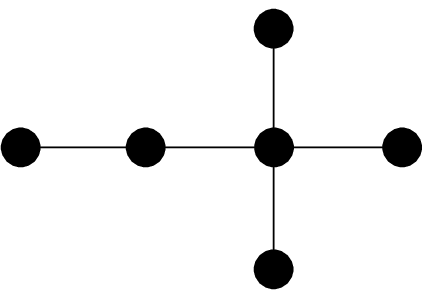,width=0.065\linewidth}}&
$\vphantom{\int\limits_a^a}$\epsfig{file=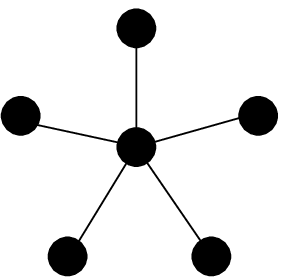,width=0.055\linewidth}$\vphantom{\int\limits_a^a}$&
\raisebox{4pt}{20} 
\\
\raisebox{3.5pt}{
\epsfig{file=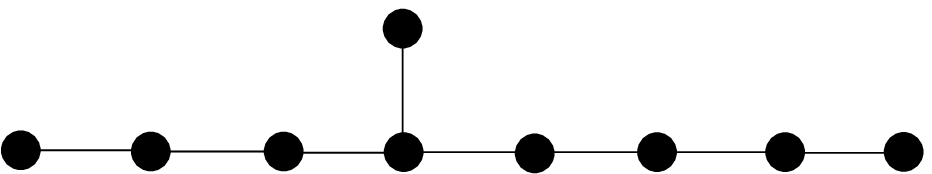,width=0.14\linewidth}}\ &
\ 
\epsfig{file=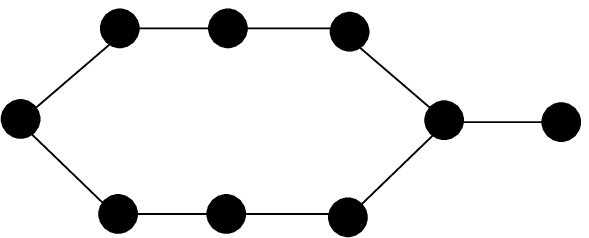,width=0.14\linewidth}\ &
\raisebox{6pt}{272} 
\\
\ 
\epsfig{file=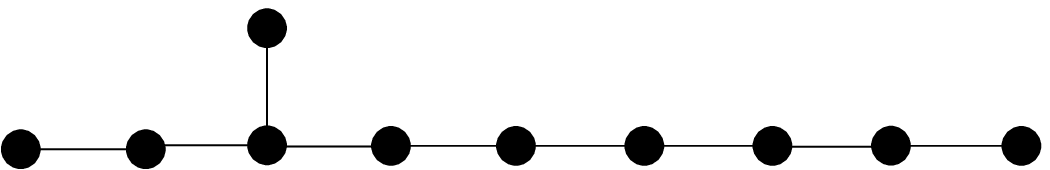,width=0.14\linewidth}\ &
\ $\vphantom{\int\limits^a_A}$
\epsfig{file=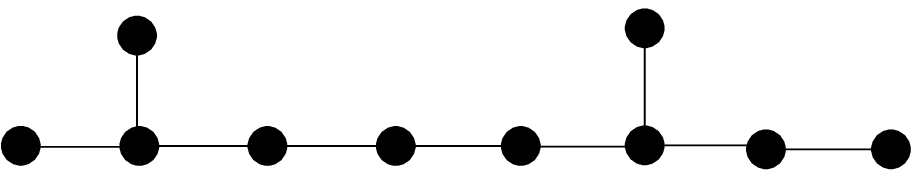,width=0.14\linewidth}$\vphantom{\int\limits^a_A}$\ &
\raisebox{1pt}{527} 
\\
\hline
\end{tabular}

\end{center}
\caption{Minimal decompositions} 
\vspace{-10pt}
\begin{center}
without decomposed dihedral angles.
\end{center}
\label{fund}
\end{table}

\noindent
{\bf Remark.\ }
The proof of the second implication does not need the
decomposition to be minimal. It will be convenient for the study of
non-minimal decompositions.

Some simplices correspond to several root systems.
Indeed, suppose that $F$ has at least one dihedral angle different
from $\frac{\pi}{2}$ and $\frac{\pi}{3}$.
Then there exists at least two ways to define the lengths of roots 
(see section~\ref{pic}). 
We will prove that for any minimal decomposition
$(F,{F_1})$ we can find a root system $\Delta$ (with fundamental
simplex $F$) such that the roots correspondent to ${F_1}$
generate a proper sublattice of $L$.
By Lemma~\ref{equiv} the condition ($*$) holds for the root system
correspondent to ${F_1}$.

First, suppose that ${F_1}$ contains exactly two copies of $F$.

\begin{lemma}
\label{2}

Suppose that $\left[W_F:W_{F_1}\right]=2$. Let $\Delta$ be any root
system with fundamental simplex $F$.
Then the  roots of $\Delta$ (or $\Delta^{\vee}$) vanishing on 
the facets of $F_1$ generate a proper sublattice $L_1$ of
$L$ (or $L^{\vee}$ respectively).
The index of the sublattice equals two, three or four.

\end{lemma}

\begin{proof}

The simplex ${F_1}$ is a union of $F$ and $F'$, where $F'$ is an image
of $F$ under the reflection with respect to some facet of $F$. 
Let $\alpha_1$ be a root vanishing on this facet. 
All but one facets of ${F_1}$ are facets of $F$.
Thus, exactly one of
$\alpha_2,...,\alpha_{n+1}$  is not orthogonal to $\alpha_1$.

Suppose that $\alpha_2$ is not orthogonal to $\alpha_1$.
We can assume that $|a_{21}|\le |a_{12}|$ 
(if $|a_{21}|\ge |a_{12}|$, consider the matrix $A^t$ instead of $A$). 
Then the facets of ${F_1}$ correspond to the roots
$\alpha_2-a_{12}\alpha_1,\alpha_2,\alpha_3,...,\alpha_{n+1}$.

Since $A$ is a matrix of hyperbolic type,
$a_{12}= -1$, $-2$, $-3$ or $-4$ ($-4$ occurs only if
${F_1}$ and ${F}$ are non-compact triangles).
In case of $a_{12}=-1$ the angle between the facets correspondent
to $\alpha_2$ and $\alpha_2-a_{12}\alpha_1$ equals $\frac{2\pi}{3}$,
and the group $W_{F_1}$ coincides with $W_F$.

Hence, we have either $a_{12}=-2$, or  $a_{12}=-3$, or  $a_{12}=-4$.
Therefore,
$L_1=\Z(\alpha_2-a_{12}\alpha_1)+\Z\alpha_2+\Z\alpha_3+\dots+
\Z\alpha_{n+1}$ is a proper sublattice of index $|a_{12}|$ of the lattice
$L=\sum\limits_{i=1}^{n+1} \Z\alpha_i$.

\end{proof}

Consider the general case.

\begin{theorem}
\label{max}
Let $W_{F_1}$ be a maximal subgroup of $W_F$.
Consider any root system $\Delta$ with fundamental simplex $F$.
Then the  roots of $\Delta$ (or $\Delta^{\vee}$) vanishing on 
the facets of $F_1$ generate a proper sublattice $L_1$ of
$L$ (or $L^{\vee}$ respectively).

\end{theorem}

\begin{proof}
Since $W_{F_1}$ is a subgroup of $W_F$,
${F_1}$ is a union of several copies of $F$,
where two copies having a common facet are symmetric with respect to
this facet.

All the minimal decompositions of simplices are described in
~\cite{Fel4-9},~\cite{Fel3} and~\cite{Klim}.
Considering these decompositions case by case,
one can find that for any minimal decomposition
$(F,{F_1})$ the simplex ${F_1}$ has at least one vertex $A$ whose
stabilizer in  $W_F$ coincides with its stabilizer in $W_{F_1}$.
In other words, all but one facets of $F_1$ are the facets of $F$
(the intersection of these facets is the vertex $A$).

Let $\alpha_1,...,\alpha_n$ be the roots vanishing on the common
facets of $F$ and ${F_1}$. Let $\alpha_{n+1}$ É $\alpha_{{F_1}}$ 
be the roots correspondent to the rest facets of $F$ and ${F_1}$ 
respectively.
The root
$\alpha_{{F_1}}$ is a linear combination 
of roots $\alpha_1,...,\alpha_{n+1}$.
The index of the sublattice generated by
$\alpha_1,...,\alpha_n,\alpha_{F_1}$ is equal to the coefficient of
$\alpha_{n+1}$.

Now take a minimal decomposition $(F,F_1)$ and a root system
$\Delta$ with a fundamental simplex $F$. Compute the coefficient of 
$\alpha_{n+1}$ when $\alpha_{{F_1}}$ is represented as a linear
combination of 
$\alpha_1,...,\alpha_{n+1}$. If the coefficient is not equal to one
then we obtain a proper sublattice (the coefficient can not be
negative, since  $\alpha_{F_1}$
is a positive root; moreover, it can not be equal to zero, otherwise simple
roots of the root system would be linearly dependent).

Suppose that  the coefficient equals one. Turn over all the arrows 
in the Dynkin diagram of $\Delta$. In other words, 
consider a root system $\Delta^{\vee}$.
A direct calculation shows that in this root system
 the coefficient we are interested in is not equal to one,
and we have a proper sublattice of $L^{\vee}$.

The sublattice is usually of index two. 
More precisely, sublattices of index different from two
occur only in the dimensions two and three
(see Fig.~\ref{2_1}--\ref{3_3}). 

\end{proof}

\section{Non-maximal subgroups}
\label{sublattice}

In this section, we prove that for any non-minimal decomposition
$(F,{F_1})$ there exist a  root system $\Delta$ 
whose simple roots vanish on the facets of $F$ and 
the root system $\Delta_1$ whose simple roots vanish on the facets of $F_1$
such that $\Delta_1$
is a root subsystem of $\Delta$. 
In general, for some decompositions $(F,{F_1})$ there exist more than
one pair of root systems
$\Delta_1\subset\Delta$ satisfying the condition described above.

\begin{lemma}\label{tower}
Let $\Delta_1$ be a root subsystem of $\Delta$,
and $\Delta_2$ be a root subsystem of $\Delta_1$. 
Then $\Delta_2$ is a root subsystem of  $\Delta$.

\end{lemma}

\begin{proof}
Suppose that
$\alpha, \beta\in\Delta_2$ and $\alpha+\beta\in\Delta$.
Since $\Delta_2\subset\Delta_1$, we have $\alpha, \beta\in\Delta_1$. 
The condition  ($*$) holds for $\Delta_1\subset\Delta$.
Thus, $\alpha+\beta\in\Delta_1$.
Since ($*$) holds for $\Delta_2\subset\Delta_1$, we have
$\alpha+\beta\in\Delta_2$.
Therefore,  ($*$) holds for $\Delta_2\subset\Delta$.

\end{proof}

The assumption of Lemma~\ref{tower} is not necessary
(see section~\ref{pic}). However, we have the following

\begin{lemma}\label{tower-}
Suppose that  $\Delta_2\subset\Delta_1\subset\Delta$, and
 $\Delta_2\subset\Delta_1$ is not a root subsystem. 
Then $\Delta_2\subset\Delta$ is not a root subsystem either.

\end{lemma}

\begin{proof}
Since  $\Delta_2\subset\Delta_1$ is not a root subsystem,
there exist
$\alpha, \beta\in\Delta_2$ such that $\alpha+\beta\in\Delta_1$ and
$\alpha+\beta\notin\Delta_2$.
Since $\Delta_1\subset\Delta$, we have $\alpha+\beta\in\Delta$. 
Therefore, ($*$) does not hold for $\Delta_2\subset\Delta$.

\end{proof}

Lemma~\ref{tower} shows  it is sufficient to
find a sequence of root systems
$\Delta_k\subset\Delta_{k-1}\subset\dots\subset\Delta_1\subset\Delta$,
such that
$\Delta_k$ corresponds to ${F_1}$, $\Delta$ corresponds to
$F$, and for any $i\le k$ the decomposition correspondent to
$\Delta_i\subset\Delta_{i-1}$ is minimal. 
Such a sequence can be constructed for almost all non-minimal decompositions.
The exclusions are two four-dimensional decompositions
(see Fig.~\ref{4_1}) and one five-dimensional decomposition
(see Fig.~\ref{5_1}). 
Root subsystems for these three decompositions are shown in
Section~\ref{pic}.

We have proved the following 

\begin{theorem}\label{real}
Let ${F}$ and ${F_1}$ be finite volume hyperbolic Coxeter simplices
having no dihedral angles different from
$\frac{\pi}{2}$, $\frac{\pi}{3}$, $\frac{\pi}{4}$, $\frac{\pi}{6}$ and $0$.
Let $W_{F_1}$ and $W_{F}$ be the groups generated by the reflections
with respect to the facets of ${F_1}$ and ${F}$ respectively.
Suppose that $W_{F_1}\subset W_{F}$.
Then there exist a root system $\Delta$ whose simple roots vanish
on the facets of $F$ and the root system $\Delta_1$ whose simple roots vanish
on the facets of $F_1$ such that $\Delta_1\subset\Delta$
is a root subsystem.

\end{theorem}

\section{Classification of maximal rank root subsystems}\label{pic}

There exist finitely many Coxeter hyperbolic simplices,
and no hyperbolic simplex exists in the dimension greater than 9
(see~\cite{29}).

Some Coxeter simplices correspond to several root systems.
To list all the root systems correspondent to Coxeter simplex, 
consider the Coxeter diagram of the
simplex and assign each multiple edge and some bold edges by an arrow
(in other words, it is sufficient to define the lengths of roots).
To obtain a Dynkin diagram of a root system,
the arrows should satisfy the only necessary condition:
if the Coxeter diagram contains a cycle without bold edges,
then  the number of arrows pointing clockwise must be equal to
the number of arrows pointing counterclockwise.
This condition should be hold by  double edges as well as by 
 triple edges (recall that an angle $\frac{\pi}{6}$
is shown in Dynkin diagram by a triple edge,
but in Coxeter diagram this angle is shown by a 4-fold edge).

If the Coxeter diagram contains a cycle with a bold edge,
the condition slightly changes (this occurs only if $n=2$).
If all the bold edges are indirected (i.e. the corresponding roots have
the same length), then the condition coincides with one 
described above. If there is an oriented bold edge, then
there are two possibilities:
either there are two bold edges with  different orientations 
and the third angle  equals  $\frac{\pi}{3}$ or 
0, or there is exactly one oriented bold edge and other two
are 2-fold edges directed to the other side.

In general case the way to list all root systems for a given Weyl group
is described in~\cite{I}.  

To obtain a complete classification of root systems  we do the
following. Consider a minimal decomposition ($F,{F_1}$).
Assign the Coxeter diagram of $F$ by arrows in all possible ways
and consider all the root systems correspondent to the simplex $F_1$.
Do this for each minimal decomposition and consider the superpositions
of minimal decompositions. This algorithm leads to the complete 
list of maximal rank hyperbolic root systems contained 
in hyperbolic root systems.

To classify regular subalgebras it is sufficient to check the condition
($*$) for each pair $\Delta_1\subset\Delta$.
In case of minimal decomposition we can use Lemma~\ref{equiv}:
it is sufficient to show that $\Delta_1$ generates a proper sublattice
of $L$.
As it was mentioned above, in case of non-minimal decompositions the
positive answer usually can be obtained applying Lemma~\ref{tower}.
Note that in this case $\Delta_1=\Delta\cap L_1$, where
$L_1$ is a root lattice for $\Delta_1$. 
Lemma~\ref{tower-} helps to make calculations shorter:
it shows immediately that ($*$) does not hold for a long list of pairs
$\Delta_1\subset\Delta$. 
In the rest cases we check ($*$) directly. There are 19 pairs of
roots systems satisfying ($*$) and corresponding to non-minimal
decompositions. It turns out that $\Delta_1=\Delta\cap L_1$ 
for all these 19 cases.
Combining this result with Lemma~\ref{equiv} and the remark to this lemma,
we have the following

\begin{theorem}\label{eq}
Let $\Delta_1\subset\Delta$ be two hyperbolic root systems of
the same rank. Let $L_1\subset L$ be the corresponding root lattices.
Then the following three conditions are equivalent:

{\sc (i)} $\Delta_1=\Delta\cap L_1$.

{\sc (ii)} $\Delta_1\subset\Delta$  is a root subsystem.

{\sc (iii)} Condition {\sc ($*$)} holds for $\Delta_1^{re}\subset\Delta^{re}$.

\end{theorem}

Below we list all the maximal rank hyperbolic root subsystems of
hyperbolic root systems.  

\begin{figure}[htb!]
\begin{center}
\psfrag{2}{\scriptsize $2$}
\psfrag{(3)}{\scriptsize $(3)$}
\epsfig{file=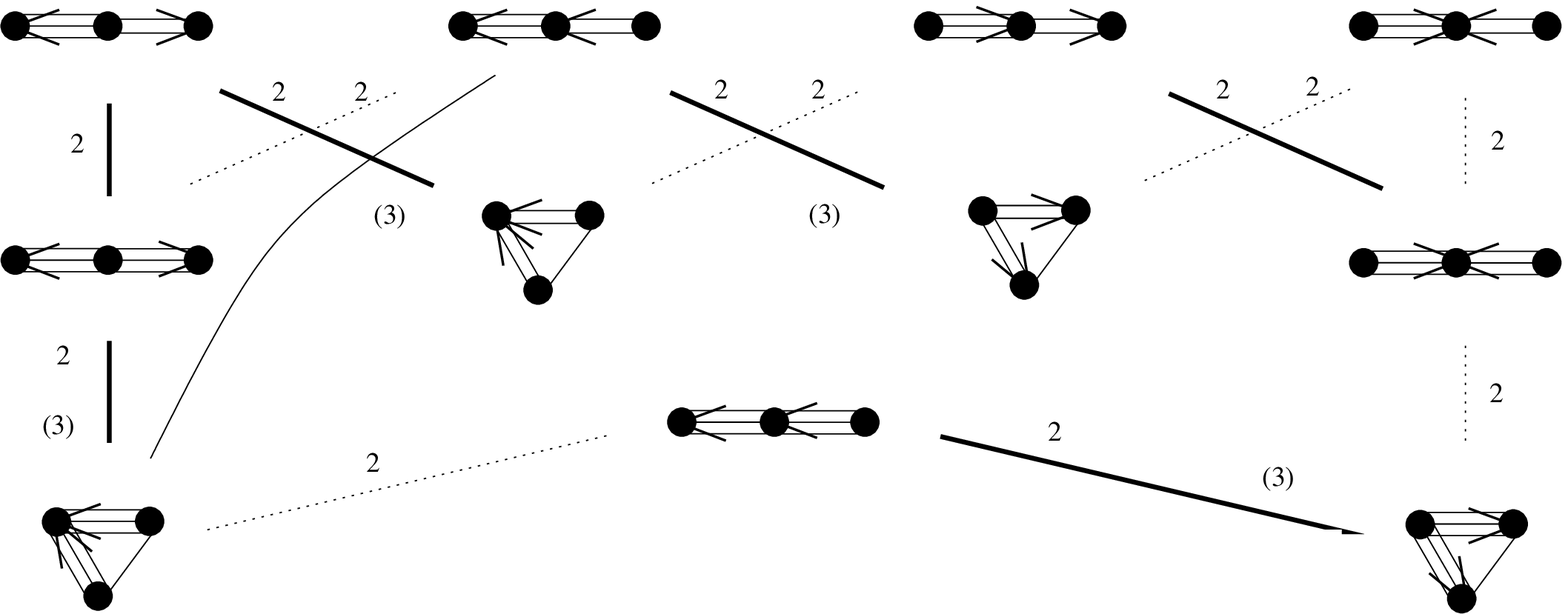,width=0.67\linewidth}\\
\end{center}
\caption{}
\label{2_1}
\end{figure}

In Fig.~\ref{2_1}--\ref{9} we use the following notation:

\noindent
Two diagrams are joined if the Weyl group correspondent to the lower
diagram is a subgroup of the Weyl group correspondent to the upper
diagram. Each edge correspondent to a minimal decomposition is
assigned with an index of the subgroup. If the decomposition is
minimal,
the lower system is a subsystem of the upper one, and the index of
the sublattice differs from 2,
then the edge is attached with the index of the sublattice (the number
in brackets).

\vspace{3pt}

\noindent
Types of edges:
\vspace{-3pt}

\begin{itemize}
\item[$$]
 \epsfig{file=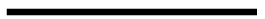,width=0.12\linewidth} \quad
decomposition is minimal, ($*$) holds;

\item[$$]
 \epsfig{file=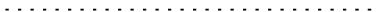,width=0.12\linewidth} \quad
decomposition is minimal, ($*$) does not hold;

\item[$$]
\epsfig{file=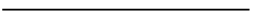,width=0.12\linewidth} \quad
\parbox[t]{9.4cm}{decomposition is non-minimal, ($*$) holds, 
but  ($*$) does not hold for at least one intermediate minimal decomposition.} 

\noindent
\end{itemize}

\noindent
A root system is not joined with a root subsystem if the
decomposition is non-minimal and  ($*$) holds for each intermediate
minimal decomposition
(see Lemma~\ref{tower}).

\vspace{-8pt}

\subsection{Triangles}

There are exactly three minimal decompositions of compact Coxeter
hyperbolic triangles having no angles different from $\frac{\pi}{2}$,
 $\frac{\pi}{3}$, $\frac{\pi}{4}$ and $\frac{\pi}{6}$.
Considering all possible lengths of simple roots, we obtain

\vspace{0pt}

\begin{figure}[htb!]
\begin{center}
\psfrag{2}{\scriptsize $2$}
\psfrag{3}{\scriptsize $3$}
\psfrag{4}{\scriptsize $4$}
\psfrag{(4)}{\scriptsize $(4)$}
\psfrag{(3)}{\scriptsize $(3)$}
\epsfig{file=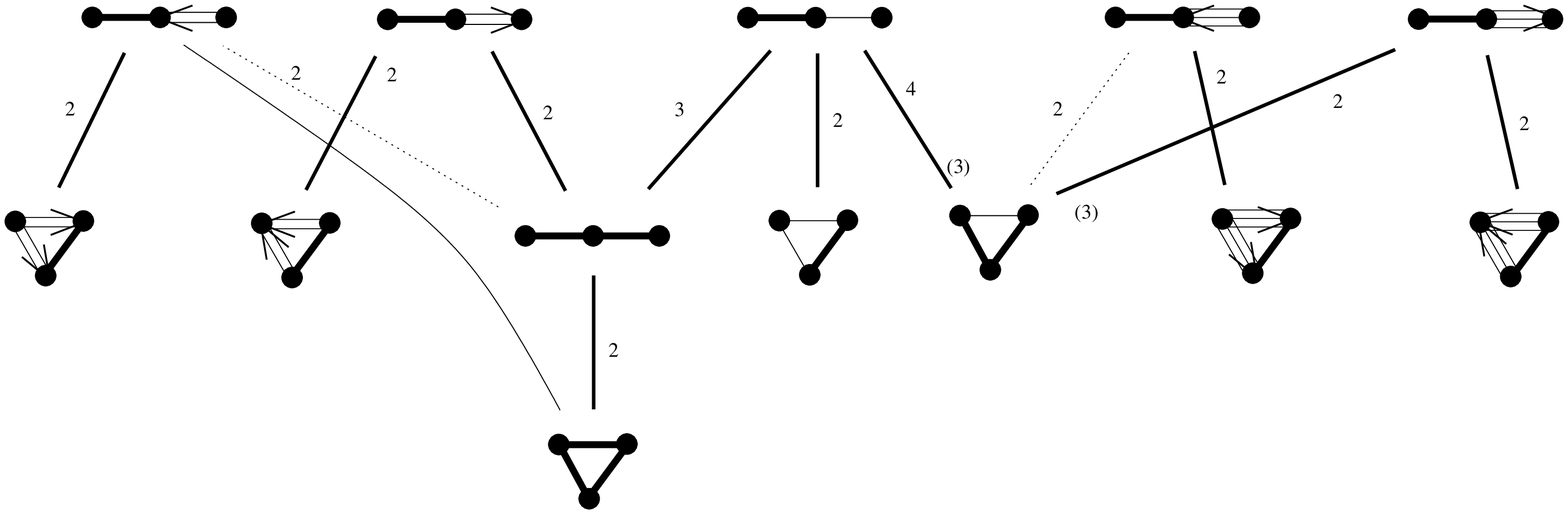,width=0.9\linewidth}\\
\end{center}
\caption{}
\label{2_2}
\end{figure}

\begin{figure}[htb!]
\begin{center}
\psfrag{2}{\scriptsize $2$}
\psfrag{3}{\scriptsize $3$}
\psfrag{4}{\scriptsize $4$}
\psfrag{(3)}{\scriptsize $(3)$}
\psfrag{(4)}{\scriptsize $(4)$}
\epsfig{file=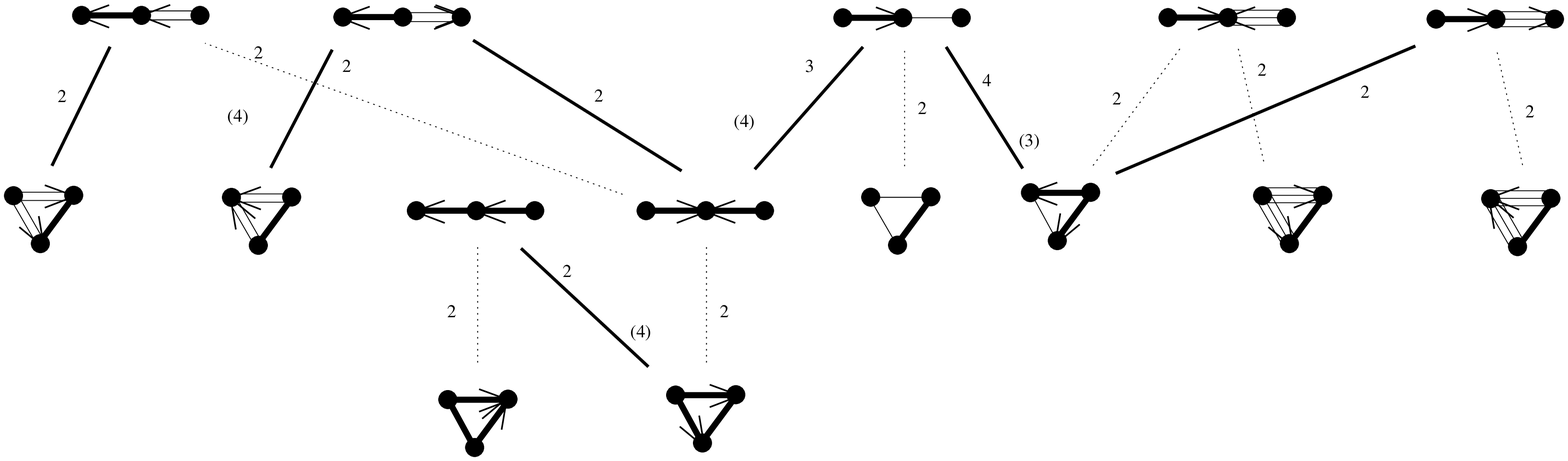,width=0.9\linewidth}\\
\end{center}
\caption{}
\label{2_3}
\end{figure}

\begin{figure}[htb!]
\begin{center}
\psfrag{2}{\scriptsize $2$}
\psfrag{3}{\scriptsize $3$}
\psfrag{4}{\scriptsize $4$}
\psfrag{(3)}{\scriptsize $(3)$}
\psfrag{(4)}{\scriptsize $(4)$}
\epsfig{file=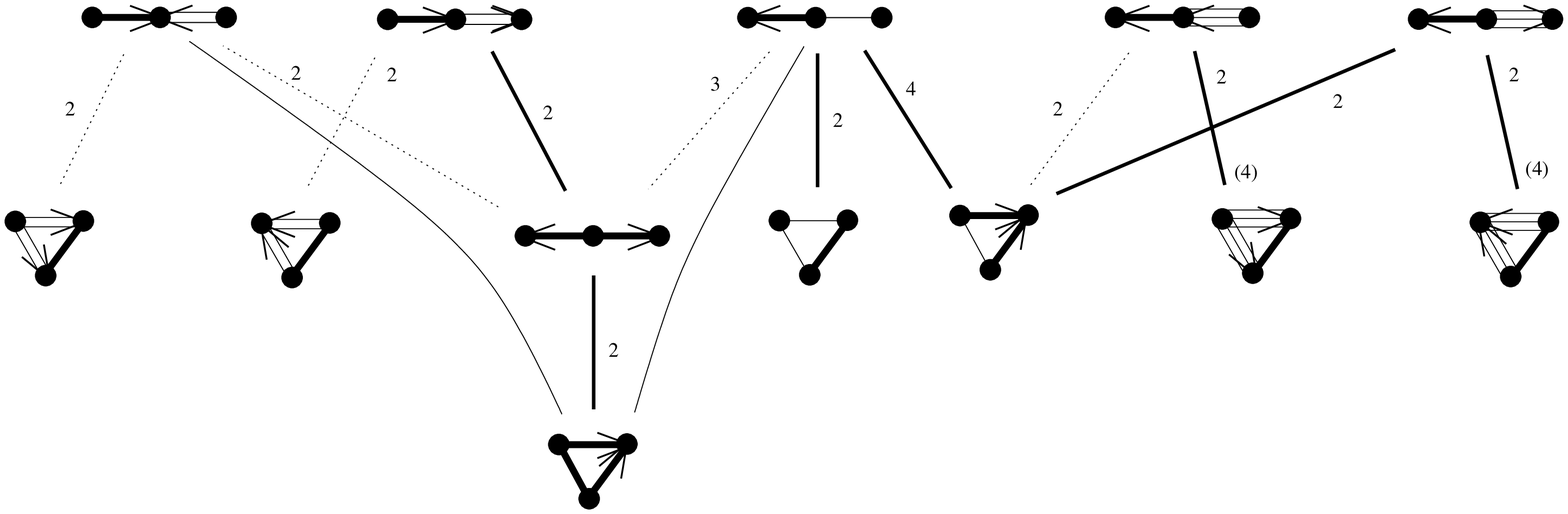,width=0.9\linewidth}\\
\end{center}
\caption{}
\label{2_4}
\end{figure}

There are six commensurability classes of non-compact Coxeter
 hyperbolic triangles
 having no  angles different from $\frac{\pi}{2}$,
 $\frac{\pi}{3}$, $\frac{\pi}{4}$, $\frac{\pi}{6}$ and $0$.
Five of these classes contain a unique triangle each,
thus, these classes produce no root subsystem.
The rest commensurability class is described in Fig.~\ref{2_2}--\ref{2_4}.

\vspace{10pt}

\subsection{Tetrahedra}

\begin{figure}[htb!]
\begin{center}
\psfrag{2}{\scriptsize $2$}
\psfrag{3}{\scriptsize $3$}
\psfrag{4}{\scriptsize $4$}
\psfrag{5}{\scriptsize $5$}
\psfrag{6}{\scriptsize $6$}
\psfrag{(3)}{\scriptsize $(3)$}
\psfrag{(4)}{\scriptsize $(4)$}
\psfrag{12}{\scriptsize $12$}
\epsfig{file=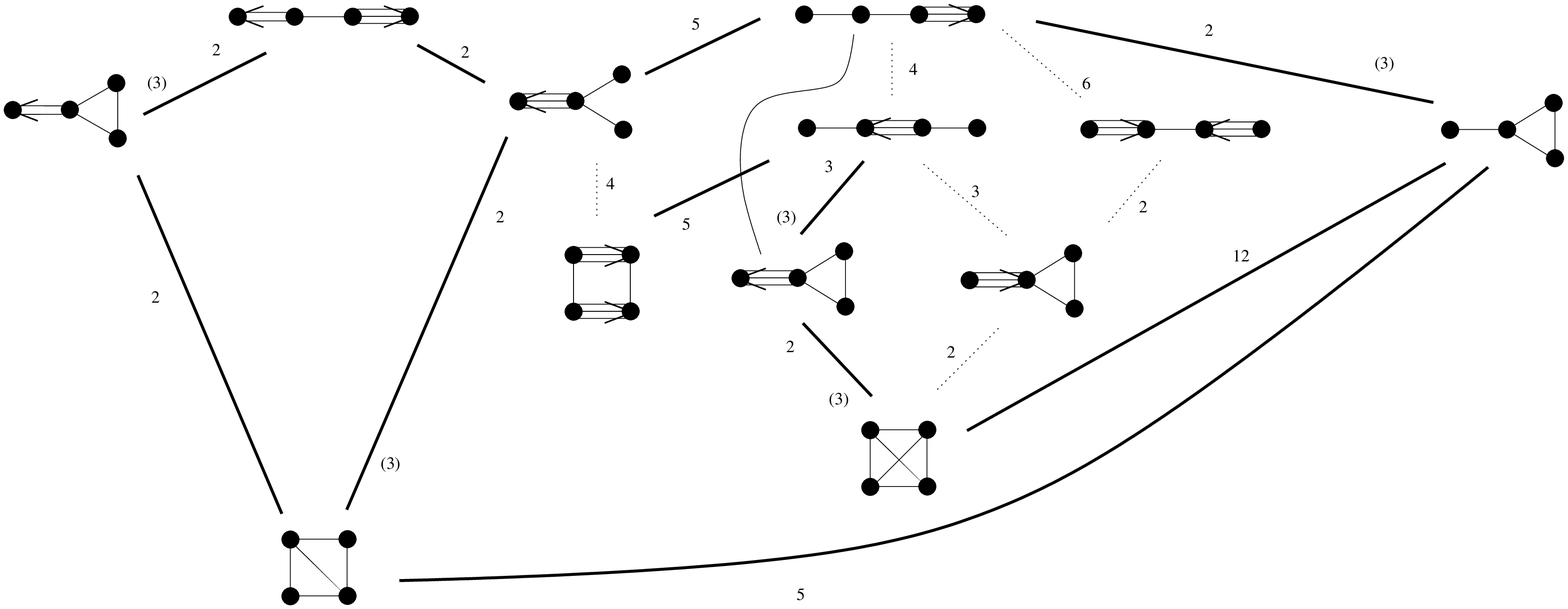,width=0.97\linewidth}\\
\end{center}
\caption{}
\label{3_1}
\end{figure}

\begin{figure}[htb!]
\begin{center}
\psfrag{2}{\scriptsize $2$}
\psfrag{3}{\scriptsize $3$}
\psfrag{4}{\scriptsize $4$}
\psfrag{5}{\scriptsize $5$}
\psfrag{6}{\scriptsize $6$}
\psfrag{12}{\scriptsize $12$}
\psfrag{(3)}{\scriptsize $(3)$}
\epsfig{file=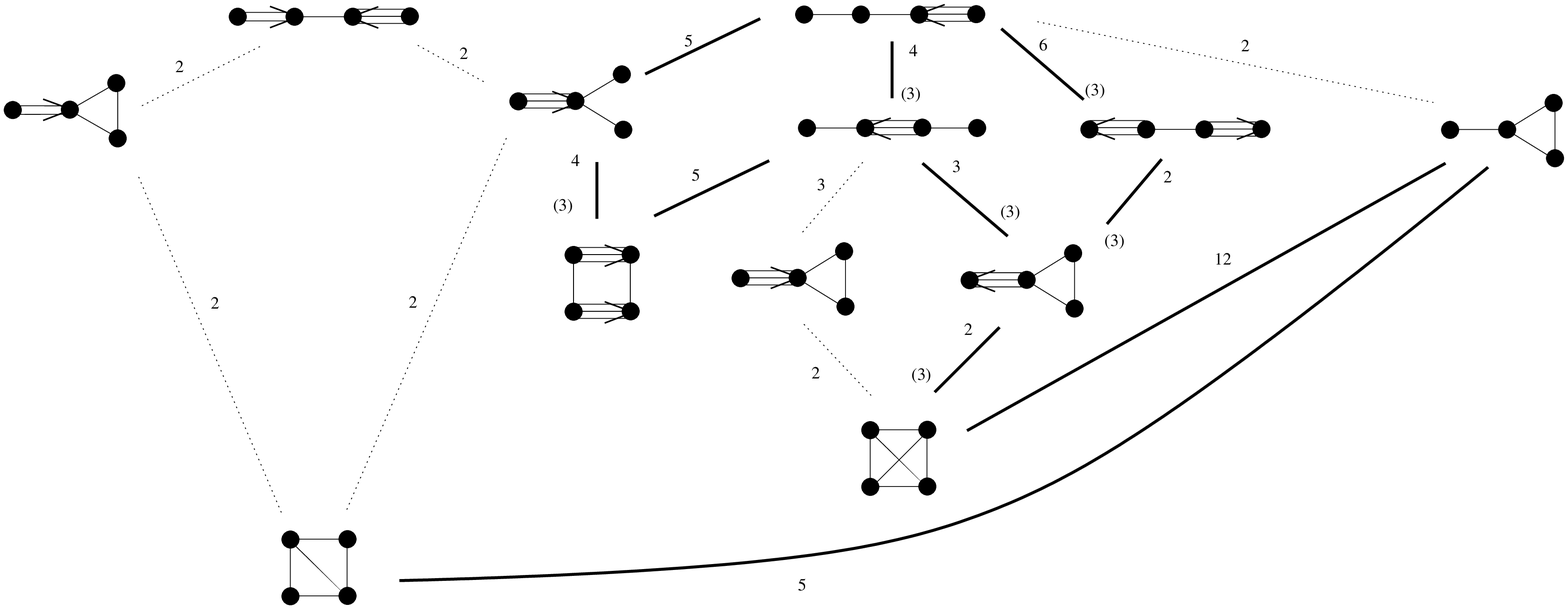,width=0.97\linewidth}\\
\end{center}
\caption{}
\label{3_2}
\end{figure}

\begin{figure}[htb!]
\begin{center}
\psfrag{2}{\scriptsize $2$}
\psfrag{3}{\scriptsize $3$}
\psfrag{4}{\scriptsize $4$}
\psfrag{5}{\scriptsize $5$}
\psfrag{6}{\scriptsize $6$}
\psfrag{12}{\scriptsize $12$}
\psfrag{(3)}{\scriptsize $(3)$}
\epsfig{file=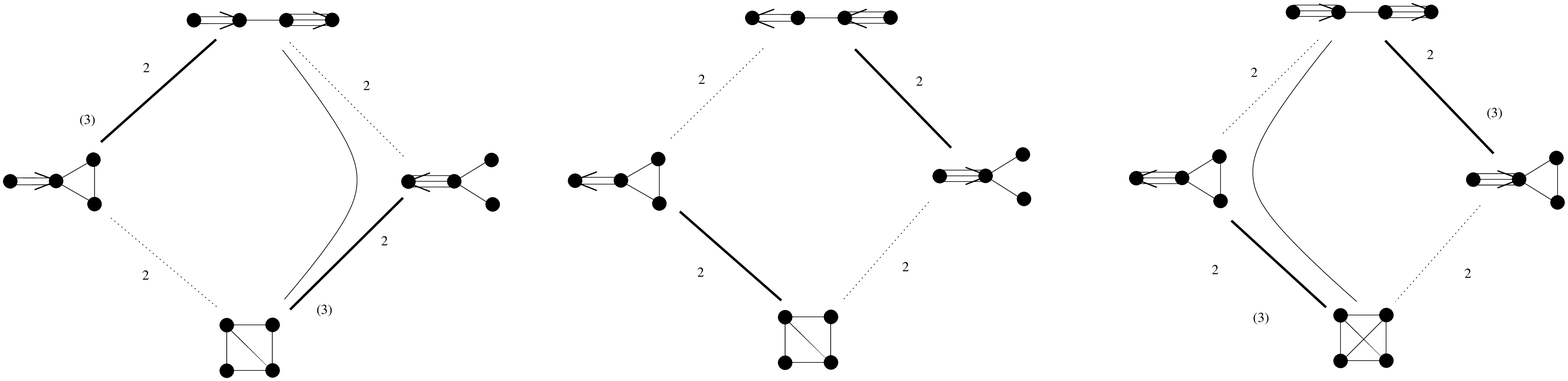,width=0.85\linewidth}\\
\end{center}
\caption{}
\label{3_3}
\end{figure}

\begin{figure}[htb!]
\begin{center}
\psfrag{2}{\scriptsize $2$}
\psfrag{3}{\scriptsize $3$}
\psfrag{4}{\scriptsize $4$}
\psfrag{5}{\scriptsize $5$}
\psfrag{6}{\scriptsize $6$}
\psfrag{(3)}{\scriptsize $(3)$}
\psfrag{(4)}{\scriptsize $(4)$}
\epsfig{file=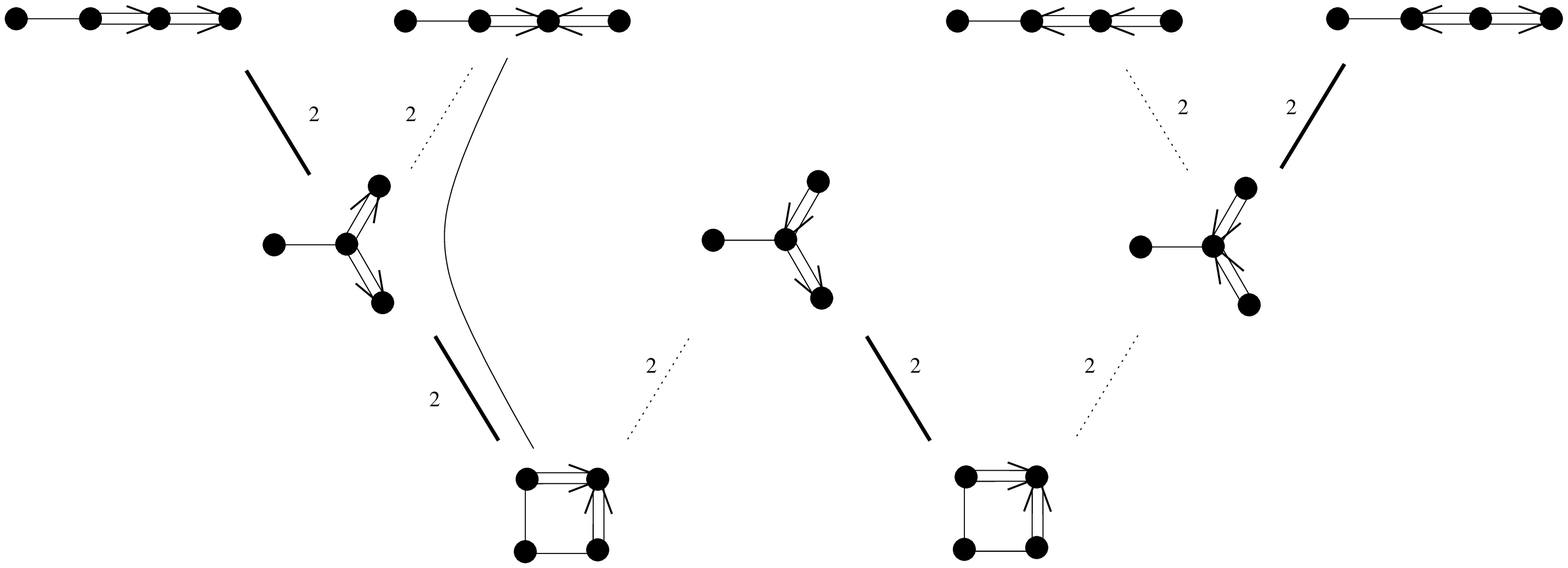,width=0.75\linewidth}\\
\end{center}
\caption{}
\label{3_5}
\end{figure}

\begin{figure}[htb!]
\begin{center}
\psfrag{2}{\scriptsize $2$}
\psfrag{3}{\scriptsize $3$}
\psfrag{4}{\scriptsize $4$}
\psfrag{5}{\scriptsize $5$}
\psfrag{6}{\scriptsize $6$}
\psfrag{(3)}{\scriptsize $(3)$}
\psfrag{(4)}{\scriptsize $(4)$}
\epsfig{file=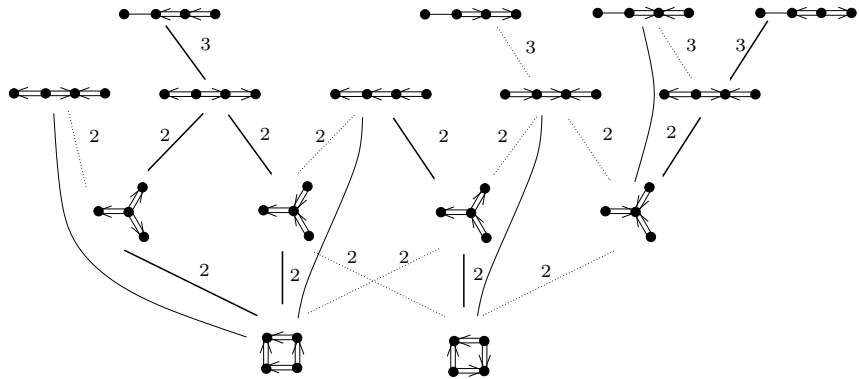,width=0.93\linewidth}\\
\end{center}
\caption{The right diagram of the second level coincides with the 
left one of the same level.}
\label{3_4}
\end{figure}

There exist four commensurability classes of  Coxeter hyperbolic tetrahedra
 having no dihedral angles different from $\frac{\pi}{2}$,
 $\frac{\pi}{3}$, $\frac{\pi}{4}$ and $\frac{\pi}{6}$
(see~\cite{Ruth}). 
Two of these classes contain a unique tetrahedron each,
thus, this classes produce no root subsystem.
The rest two commensurability classes are described in 
Fig.~\ref{3_1}--\ref{3_3} and Fig.~\ref{3_5}--\ref{3_4}.

\vspace{13pt}

\vspace{13pt}

\noindent

To obtain diagrams in Fig.~\ref{3_2} one can turn over all the arrows
on the diagrams shown in Fig.~\ref{3_1}.

\vspace{8pt}

\vspace{6pt}

The rest possibilities to assign the arrows (that correspond  
to the root systems containing real roots of three different lengths)
are shown in Fig.~\ref{3_3}.

Root systems correspondent to the second commensurability class
are shown in Fig.~\ref{3_5} and~\ref{3_4}.

\subsection{Four-dimensional simplices}

\begin{figure}[htb!]
\begin{center}
\psfrag{2}{\scriptsize $2$}
\psfrag{3}{\scriptsize $3$}
\psfrag{4}{\scriptsize $4$}
\psfrag{5}{\scriptsize $5$}
\psfrag{6}{\scriptsize $6$}
\psfrag{10}{\scriptsize $10$}
\epsfig{file=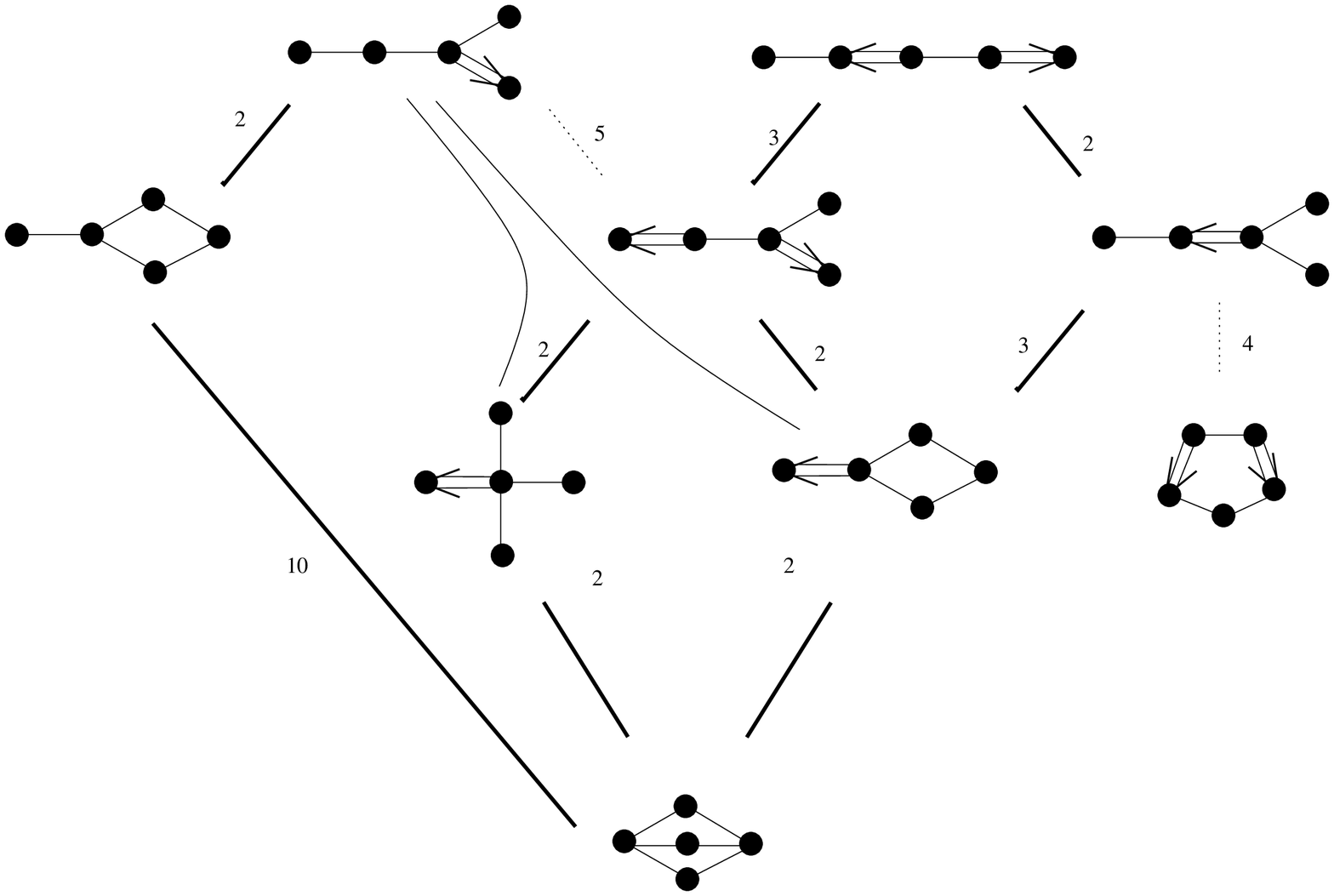,width=0.88\linewidth}\\
\end{center}
\caption{}
\label{4_1}
\end{figure}

\begin{figure}[htb!]
\begin{center}
\psfrag{2}{\scriptsize $2$}
\psfrag{3}{\scriptsize $3$}
\psfrag{4}{\scriptsize $4$}
\psfrag{5}{\scriptsize $5$}
\psfrag{6}{\scriptsize $6$}
\psfrag{10}{\scriptsize $10$}
\epsfig{file=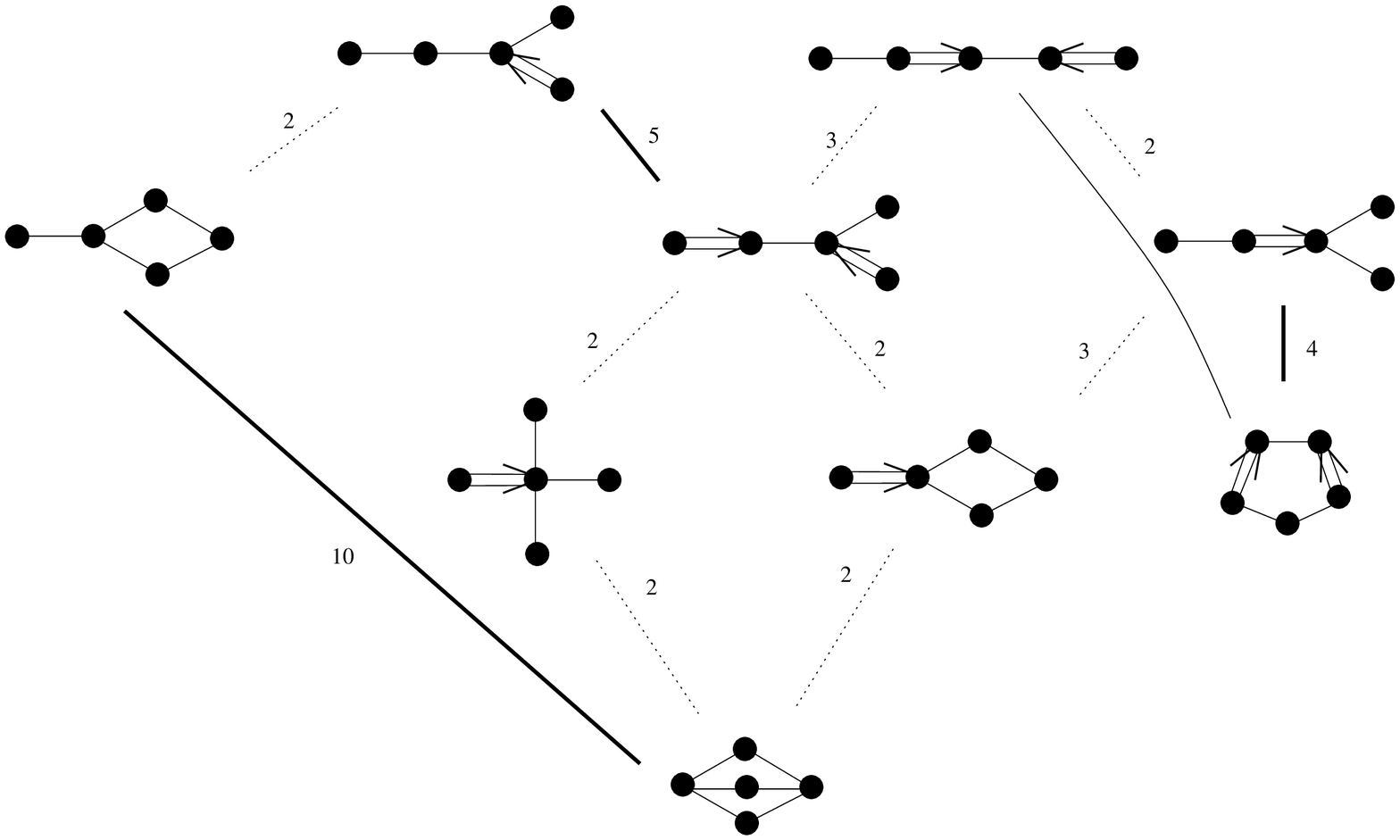,width=0.85\linewidth}\\
\end{center}
\caption{}
\label{4_2}
\end{figure}

\begin{figure}[htb!]
\begin{center}
\psfrag{2}{\scriptsize $2$}
\psfrag{3}{\scriptsize $3$}
\psfrag{4}{\scriptsize $4$}
\psfrag{5}{\scriptsize $5$}
\psfrag{6}{\scriptsize $6$}
\epsfig{file=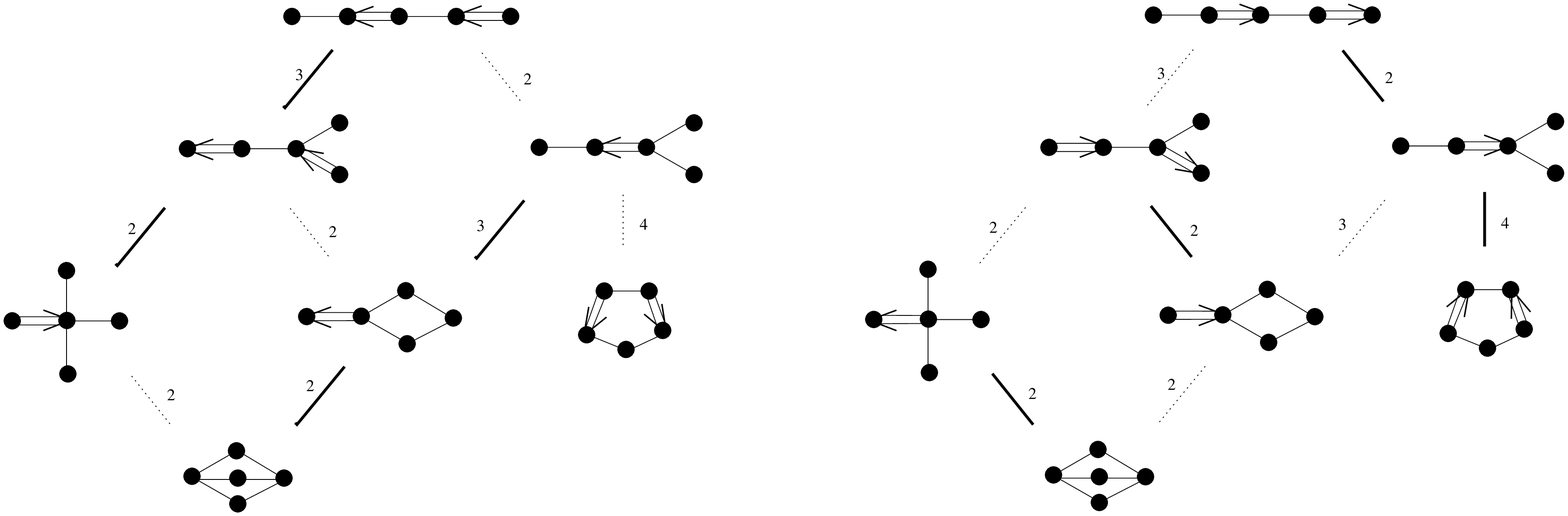,width=0.99\linewidth}\\
\end{center}
\caption{}
\label{4_3}
\end{figure}

According to~\cite{Ruth} (see also~\cite{Fel4-9}), 
there are exactly two commensurability classes of four-dimensional
Coxeter hyperbolic simplices
 having no dihedral angles different from $\frac{\pi}{2}$,
 $\frac{\pi}{3}$, $\frac{\pi}{4}$ and $\frac{\pi}{6}$. 
One of these classes contains a unique simplex,
thus, this class produces no root subsystem.
The rest commensurability class is described in 
Fig.~\ref{4_1}--\ref{4_3}.

\noindent
To obtain diagrams in Fig.~\ref{4_2} one can turn over all the arrows
on the diagrams shown in Fig.~\ref{4_1}.

The rest possibilities to assign the arrows (that correspond 
to the root systems containing real roots of three different lengths)
are shown in Fig.~\ref{4_3}.

\vspace{5pt}

\subsection{Five-dimensional simplices}

\begin{figure}[htb!]
\begin{center}
\psfrag{2}{\scriptsize $2$}
\psfrag{3}{\scriptsize $3$}
\psfrag{4}{\scriptsize $4$}
\psfrag{5}{\scriptsize $5$}
\psfrag{6}{\scriptsize $6$}
\psfrag{8}{\scriptsize $8$}
\psfrag{10}{\scriptsize $10$}
\psfrag{20}{\scriptsize $20$}
\epsfig{file=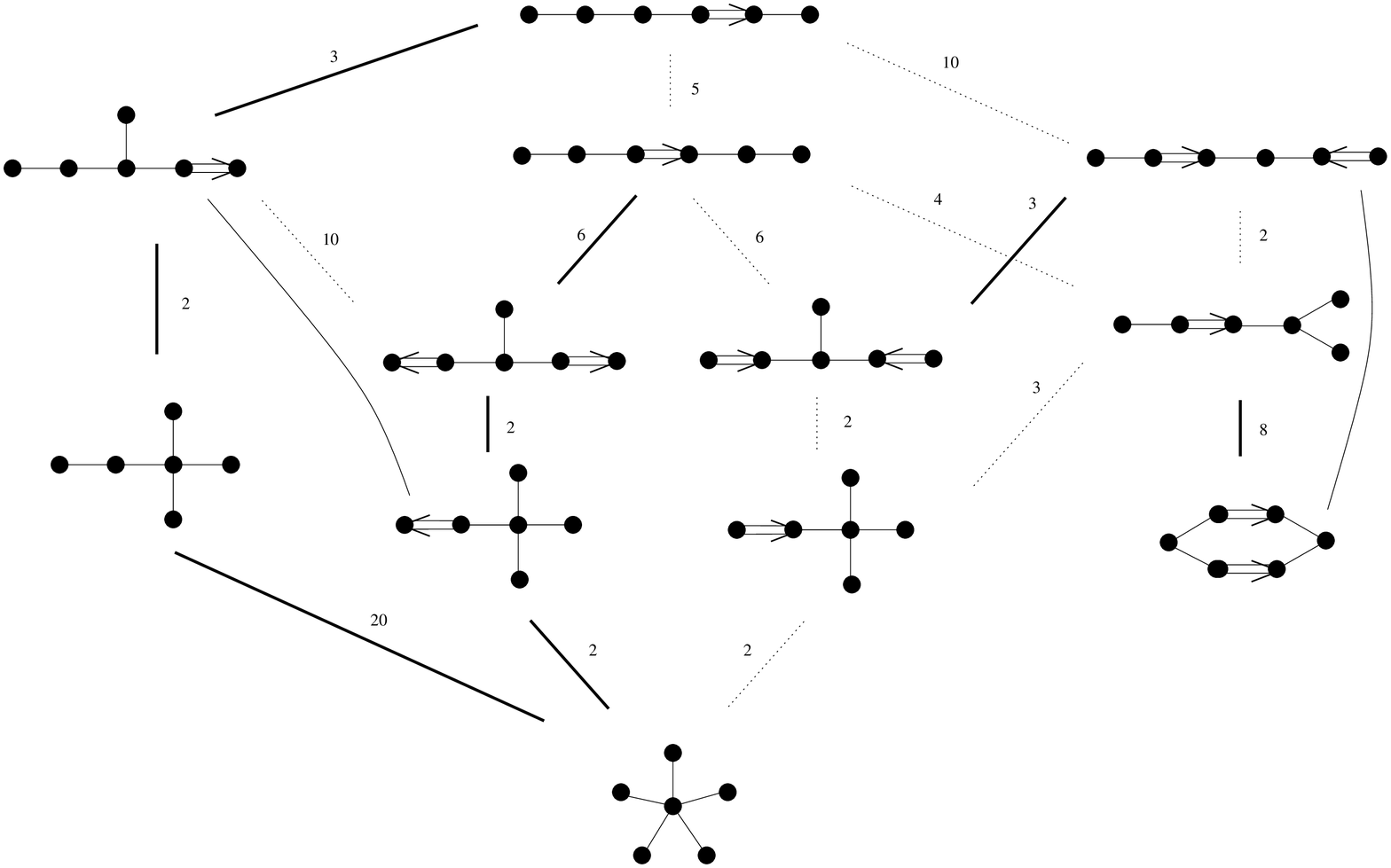,width=0.92\linewidth}\\
\end{center}
\caption{}
\label{5_1}
\end{figure}

\begin{figure}[htb!]
\begin{center}
\psfrag{2}{\scriptsize $2$}
\psfrag{3}{\scriptsize $3$}
\psfrag{4}{\scriptsize $4$}
\psfrag{5}{\scriptsize $5$}
\psfrag{6}{\scriptsize $6$}
\psfrag{8}{\scriptsize $8$}
\psfrag{10}{\scriptsize $10$}
\psfrag{20}{\scriptsize $20$}
\epsfig{file=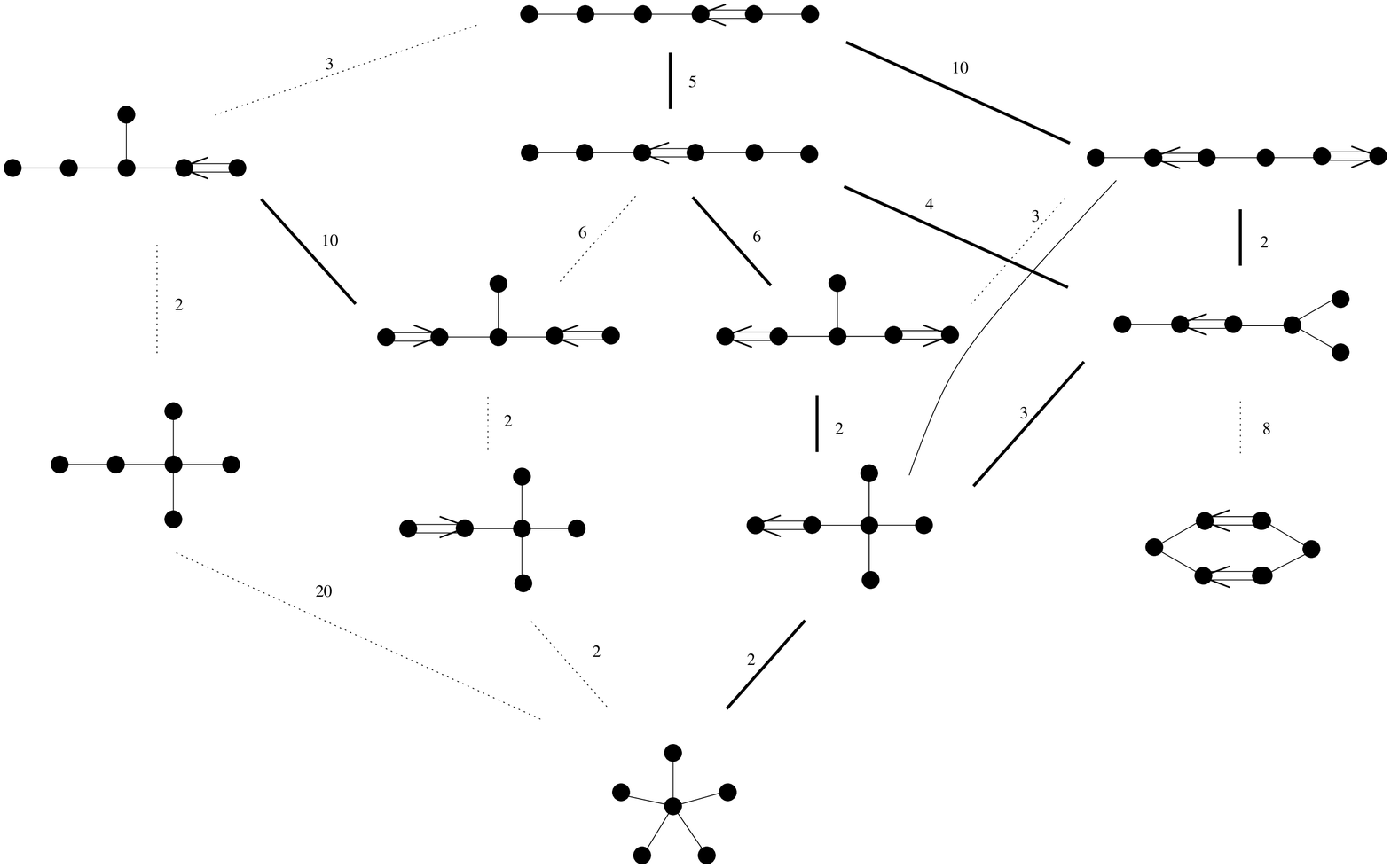,width=0.87\linewidth}\\
\end{center}
\caption{}
\label{5_2}
\end{figure}

\begin{figure}[htb!]
\begin{center}
\psfrag{2}{\scriptsize $2$}
\psfrag{3}{\scriptsize $3$}
\psfrag{4}{\scriptsize $4$}
\psfrag{5}{\scriptsize $5$}
\psfrag{6}{\scriptsize $6$}
\psfrag{8}{\scriptsize $8$}
\psfrag{10}{\scriptsize $10$}
\psfrag{20}{\scriptsize $20$}
\epsfig{file=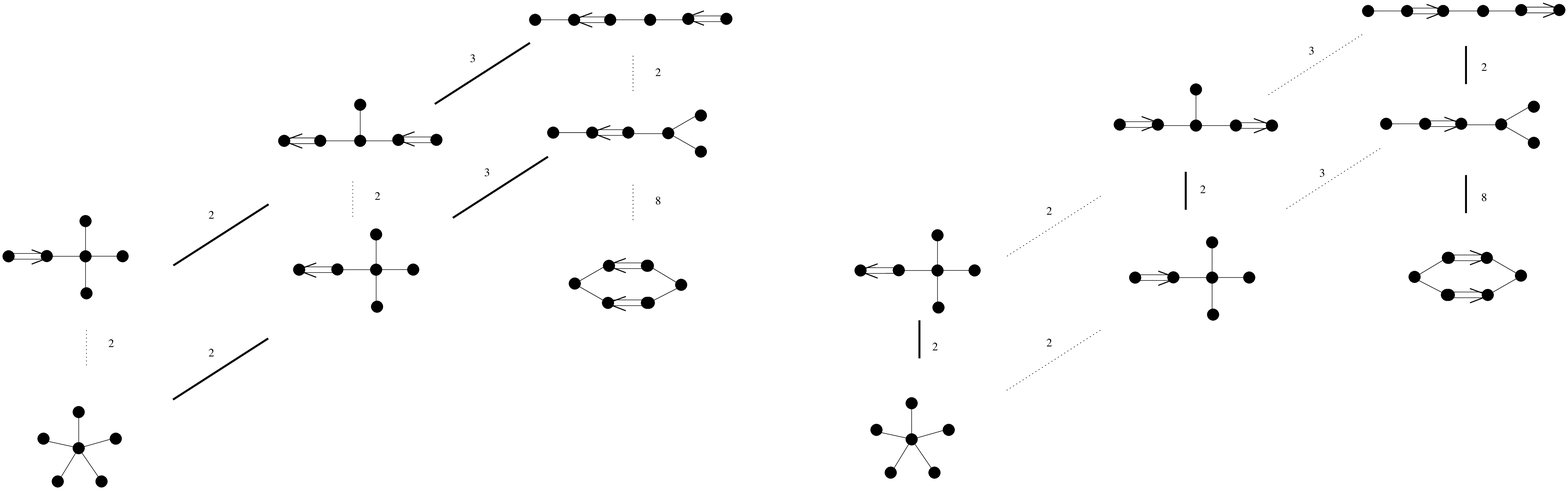,width=0.96\linewidth}\\
\end{center}
\caption{}
\label{5_3}
\end{figure}

According to~\cite{Ruth} (see also~\cite{Fel4-9}), 
there are exactly three commensurability classes of five-dimensional
Coxeter hyperbolic simplices.
Two of these classes contain a unique simplex each,
thus, these classes produce no root subsystem.
The rest commensurability class is described in 
Fig.~\ref{5_1}--\ref{5_3}.

\vspace{15pt}

\vspace{2pt}

To obtain diagrams in Fig.~\ref{5_2} one can turn over all the arrows
on the diagrams shown in Fig.~\ref{5_1}.

\vspace{2pt}

\vspace{4pt}

The rest possibilities to assign the arrows (that correspond 
to the root systems containing real roots of three different lengths)
are shown in Fig.~\ref{5_3}.

\subsection{Simplices of dimensions 6-9}

\begin{figure}[htb!]
\begin{center}
\psfrag{2}{\scriptsize $2$}
\epsfig{file=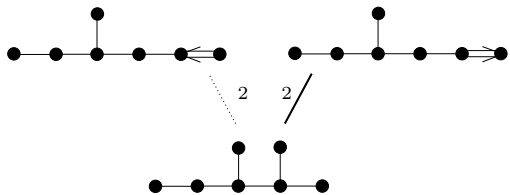,width=0.55\linewidth}\\
\end{center}
\caption{Six-dimensional root subsystems.}
\label{6}
\end{figure}

\begin{figure}[htb!]
\begin{center}
\psfrag{2}{\scriptsize $2$}
\epsfig{file=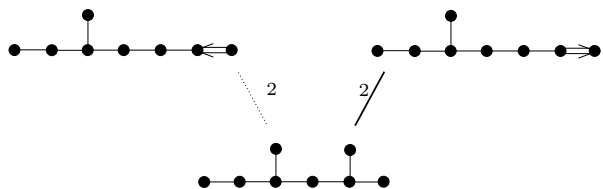,width=0.65\linewidth}\\
\end{center}
\caption{Seven-dimensional root subsystems.}
\label{7}
\end{figure}

\begin{figure}[htb!]
\begin{center}
\psfrag{2}{\scriptsize $2$}
\psfrag{272}{\scriptsize $272$}
\epsfig{file=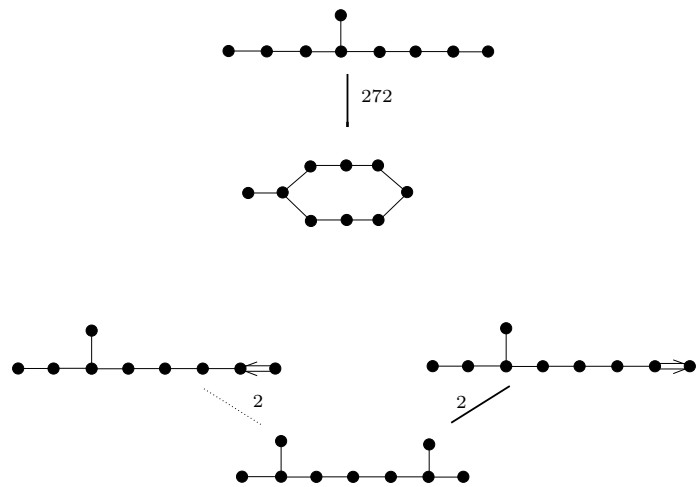,width=0.75\linewidth}\\
\end{center}
\caption{Eight-dimensional root subsystems.}
\label{8}
\end{figure}

\begin{figure}[htb!]
\begin{center}
\psfrag{2}{\scriptsize $2$}
\psfrag{527}{\scriptsize $527$}
\epsfig{file=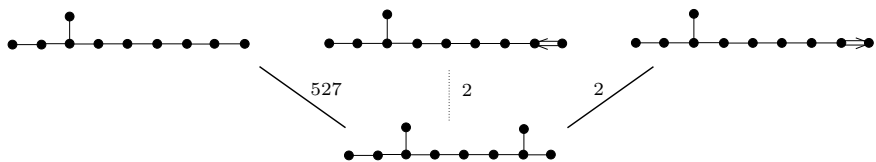,width=0.95\linewidth}\\
\end{center}
\caption{Nine-dimensional root subsystems.}
\label{9}
\end{figure}

There are no non-minimal decompositions in the 
dimensions 6-9 (see~\cite{Ruth}).
All root subsystems in these dimensions are shown in Fig.~\ref{6}--\ref{9}.


\pagebreak

\vspace{25pt}

\noindent
Independent Univ. of Moscow,\\
e-mail:\qquad pasha@mccme.ru

\end{document}